\documentclass{article}[12pt]

\usepackage[T2A]{fontenc}
\usepackage[cp1251]{inputenc}
\usepackage[english]{babel}
\usepackage{amsmath}
\usepackage{amssymb}
\usepackage{amsthm}
\usepackage{mathrsfs}

\usepackage{hyperref}

\theoremstyle{en-defi}
\newtheorem{ru-defi}{Определение}
\theoremstyle{plain}
\newtheorem{ru-remark}{Замечание}
\newtheorem{ru-theorem}{Теорема}
\newtheorem{ru-proposition}{Предложение}

\newtheorem{ru-corollary}{Следствие}
\newtheorem{ru-lemma}{Лемма}

\newlength{\wdth}

\newtheorem{en-defi}{Definition}
\theoremstyle{plain}
\newtheorem{en-remark}{Remark}
\newtheorem{en-lemma}{Lemma}
\newtheorem{en-theorem}{Theorem}
\newtheorem{en-proposition}{Proposition}
\newtheorem{en-corollary}{Corollary}
\newtheorem{en-OldTheorem}{Theorem}

\newtheorem{en-Theorem}{Theorem}

\newif\ifen
\entrue 

\newif\ifru
\rufalse

\begin{document}

\ifen
\title{On conditions for Dobrushin's Central limit theorem for non-homogeneous Markov chains}
\author{A.\,Yu.~Veretennikov\footnote{A.A. Kharkevich Institute for Information Transmission Problems of Russian Academy of Sciences,  Moscow, Russian Federation \& S.M. Nikol'skii Mathematical Institute, RUDN University, Moscow, Russian Federation}, A.\,I.~Nurieva\footnote{National Research University Higher School of Economics, Moscow, Russian Federation \& A.A. Kharkevich Institute for Information Transmission Problems of Russian Academy of Sciences, Moscow, Russian Federation}}





\fi

\maketitle

\ifen
\abstract{A new sufficient condition is proposed in the problem of Central Limit Theorem in the array scheme for non-homogeneous Markov Chains, with a possibility that the minimum of Markov--Dobrushin's ergodic coefficient may be closer to zero than it is required for Dobrushin's main condition.}

\medskip


\medskip

\fi

\makeatletter

\ifen
\section{Introduction}
A new version of Dobrushin's central limit theorem (CLT in what follows) for non-homogeneous Markov chains is proposed in this paper. Prior to the seminal article \cite{Dobr56}, earlier results in this area were established in \cite{Markov, Bern,  Sapogov, Linnik}. R.L. Dobrushin found that the following sufficient condition for the CLT in the array scheme is very close to being necessary (under the additional non-degeneracy variances condition, see what follows), 
$$
\lim_{n \to \infty} n^{1/3} \alpha_n = \infty,
$$
where $\alpha_n$ is Markov -- Dobrushin's (MD) ergodic coefficient (see (\ref{dpi}) for the definition in the next section). The optimality of this value $1/3$ was confirmed by the Bernstein -- Dobrushin example (see \cite[ch. 2]{Bern}. Among various extensions of this result,  see the short communication \cite{Shi} for Markov chains with memory. In \cite{Seth-Varadhan} (2005) Dobrushin's results were re-proved by a much simpler martingale method. This paper uses the same technique as in \cite{Seth-Varadhan} with certain modifications, and establishes that, in fact,  $\alpha_n$ may be much closer to zero than it is required for Dobrushin's main condition, but the CLT still holds under certain little additional changes in the assumptions.  
The motivation of a CLT does not require any explanation. Some more recent references on the theme of Dobrushin's CLT may be found here: \cite{Peligrad12,SinnChen13,ArlottoSteele16, Hafouta23}.

The paper consists of 5 sections: Introduction, The setting and earlier results, Main results, Auxiliary results, Proof of the theorem.

\fi

\ifen
\section{The setting and earlier results}

Let us consider a non-homogeneous Markov chain on a probability space $(\Omega, {\cal F}, \mathsf P)$. Let $\pi = \pi_{i,i+1} (x, dy)$ be the Markovian transition kernel on $(\mathbf{X}, \mathcal{B}(\mathbf{X}))$ at time $i$, where $(\mathbf{X}, \mathcal{B}(\mathbf{X}))$ is a measurable space with standard requirements on the Markov process (see \cite{Dynkin}). The transition probabilities are then given by the formula
$
\mathsf P[X_{i+1} \in A | X_i = x]\stackrel{\text{a.s.}} = \pi_{i,i+1}(x,A),
$
where $\{X_i: 1 \leq i \leq n \}$ is the trajectory of the process on  $[1,n]$. In particular, if the initial distribution is $X_1 \sim \mu$ (it is convenient here to regard $1$ as the initial time, not $0$) then the distribution at time $k \geq 2$ is given by the measure $\mu\pi_{1,2}\pi_{2,3} \dots \pi_{k-1,k}$.

For $i < j = i+k$ let us denote by $\pi_{i,i+k}$ the $k$-step transition kernel, which has a form due to the Chapman-Kolmogorov equation,
$$
\pi_{i,j} = \pi_{i,i+1}\pi_{i+1,i+2} \dots \pi_{j-1,j}.
$$

\begin{en-defi}\label{contr_coef}
Markov -- Dobrushin's coefficient (the MD coefficient in what follows) $\delta(\pi)$ is defined by any of the three equivalent expressions, 
\begin{equation}\label{dpi}
\begin{aligned}
&\!\delta(\pi) 
\! = \!\sup_{\substack{x_1, x_2 \in \mathbf{X},\\ A \in \mathcal{B}(\mathbf{X})} }|\pi(x_1, A) \!-\! \pi(x_2, A)| 
\!=\! \frac{1}{2} \sup_{\substack{x_1, x_2 \in \mathbf{X},\\ \|f\|_{B}\le 1}} | \int f(y) [\pi (x_1, dy) \!-\! \pi (x_2, dy)]|
 \\\\
& \!= \sup_{\substack{x_1,x_2 \in \mathbf{X}, 
 \\ 
u \in \mathcal{U}}} |(\pi u)(x_1) \!-\! (\pi u)(x_2)|, 
\quad \text{where} \quad \mathcal{U} \!=\! \{ u(\cdot): \sup_{y_1, y_2}|u(y_1) \!-\! u(y_2)| \leq 1\}.
\end{aligned}
\end{equation}
\end{en-defi}
\noindent
Here $\|f\|_{B} = \sup_{x\in {\mathbf X}}|f(x)|$.
It follows from the definition that $0 \leq \delta(\pi) \leq 1$, and that $\delta(\pi) = 0$ if and only if the measure $\pi (x, dy)$ does not depend on $x$, and in the latter case the second random variable does not depend on the first one. 
To denote integrals, standard notations will be used, 
$\displaystyle (\mu \pi)(A) = \int \mu(dx) \pi(x, A)$ and $\displaystyle (\pi u)(x) = \int  u(y) \pi(x, dy)$. Let  
$$
\operatorname{Osc}(u) := \sup_{x_1, x_2}|u(x_1) - u(x_2)|.
$$ 

In what follows, this semi-norm applied to random variables (see lemmata \ref{lem3}, \ref{lem4}) will be understood in the sense 
\begin{align}\label{Osc-rv}
&\operatorname{Osc}(X) := 
\mathop{\text{ess}}\limits_{\mathsf P} 
\sup_{\omega\in \Omega} X(\omega) - \mathop{\text{ess}}\limits_{\mathsf P} 
\inf_{\omega\in \Omega} X(\omega).
\end{align}

Further, for any transition kernels $\pi_1, \pi_2$ the following inequality holds true, 
$
\delta(\pi_1 \pi_2) \leq \delta(\pi_1) \delta(\pi_2)
$, 
where $\pi_1 \pi_2$ is the two-step transition kernel. 
Similarly, by induction,
\begin{equation}\label{piij}
\delta(\pi_{i,j}) \leq \delta(\pi_{i,i+1})\delta(\pi_{i+1,i+2}) \dots \delta(\pi_{j-1,j}).
\end{equation}
Accordingly, the $k$-step coefficient $\alpha(\pi_{i,i+k})$ is defined as 
\begin{equation}\label{apik}
\alpha(\pi_{i,i+k}) := 1 - \delta(\pi_{i,i+k}).
\end{equation}

\noindent
In what follows, $\mathsf{D} [Y]$ stands for the variance of a random variable $Y$.  
In the array scheme, Markov processes  $\{ X_i^{(n)}: 1 \leq i \leq n \}$ are defined for any $n \geq 1$ on a  probability space $(\Omega, {\cal F}, \mathsf P)$ with state space $\mathbf{X}$, with the transition kernels  $\{\pi^{(n)}_{i,i+1} = \pi^{(n)}_{i,i+1} (x, dy): 1 \leq i \leq n-1 \}$ and with the initial distribution $\mu^{(n)}$. Let
\begin{equation}\label{an}
\alpha_n := \min_{1 \leq i \leq n -1} \alpha(\pi^{(n)}_{i,i+1}).
\end{equation}
Further, let  $\{ f^{(n)}_i: \mathbf{X} \mapsto \mathbb R, \, 1 \leq i \leq n \}$ be real-valued functions on $\mathbf{X}$, and for any $n \geq 1$ let
$$
S_n := \sum_{i=1}^n  f^{(n)}_i(X^{(n)}_i).
$$
The main result of Dobrushin's paper \cite{Dobr56} along with the main corollary states weak convergence (see also \cite{Seth-Varadhan}),
\[
\frac{S_n - \mathsf{E}(S_n)}{\sqrt{{\mathsf D}(S_n)}} \implies {\cal N}(0,1),
\]
under the conditions that (in a simplified form)
\begin{equation}\label{Dobr-cond}
\sup\limits_{1 \leq i \leq n} \sup\limits_{x \in \mathbf{X}} |f^{(n)}_i(x)| <\infty, \;  \& \; \inf_{i,n}\mathsf{D}(f^{(n)}_i(X^{(n)}_i))>0, \; \& \; \lim\limits_{n \to \infty} n \alpha_n^3 = \infty.
\end{equation}
In our main result this condition will be relaxed, allowing significantly smaller values of $\alpha_n$.

\fi

\ifen
\section{Main results}\label{sec3}
Notice that 
(\ref{piij}), (\ref{apik}) and (\ref{an}) imply 
the following elementary inequality: 
\begin{equation}\label{e2-1}
\delta(\pi_{i,j}) \leq (1 - \alpha_n)^{j-i}.
\end{equation}

Let $(\beta^{(n)}_k,\, 1\le k\le n)$ for any $n\ge 1$ be a non-random sequence of ones and zeros, like $(110010\ldots01)$, and let for $1\le i \le j\le n$,
$$
\kappa^\beta_n = \sum_{k=1}^{n} \beta^{(n)}_k, 
\;  
\kappa^\beta_{n,j} = \sum_{k=1}^{j} \beta^{(n)}_k, 
\;  
\kappa^\beta_{n,i,j} = \sum_{k=i}^{j} \beta^{(n)}_k.
$$ 
Denote
$$
\alpha_n^\beta :=  \min_{1 \leq i \leq n -1: \, \beta^{(n)}_i=1} \alpha(\pi^{(n)}_{i,i+1}).
$$
Then (\ref{e2-1}) may be slightly weakened to
\begin{equation*}
\delta(\pi_{i,j}) \leq (1 - \alpha^\beta_n)^{\kappa^\beta_{n,i,j}}.
\end{equation*}

\begin{en-defi}\label{defb}
We say that the condition $(H_\beta)$ is met if there exist $m_0>0$ and $c>0$ such that for any $n$ there exists a sequence $(\beta^{(n)}_k,\, 1\le k\le n)$ with the following property: 
for any pair $(i<j)$ with $j-i\ge m_0$, the inequality 
\begin{equation}\label{ijm0}
k_j^\beta - k_i^\beta \ge c (j-i)
\end{equation}
is satisfied.
\end{en-defi}
In what follows such a sequence for each $n$ is fixed.

\begin{en-theorem}\label{thm1}
Let the assumption $(H_\beta)$ is met, and let also
$$
\sup_{1 \leq i \leq n} \sup_{x \in \mathbf{X}} |f^{(n)}_i(x)| \leq C_n <\infty. 
$$
Then, if condition 
\begin{equation*}
\lim_{n \to \infty} C^2_n (\alpha_n)^{-1} \, (\alpha^\beta_n)^{-2} \big[ \sum_{i=1}^n  \mathsf{D}(f^{(n)}_i(X^{(n)}_i)) \big]^{-1} = 0
\end{equation*}
is satisfied, and if 
there exist $m_0>0$ and $c>0$ such that for any couple $(i<j)$ with $j-i\ge m_0$, the assumption $(H_\beta)$ holds, 
then the CLT holds true, 
\begin{equation} \label{CLT}
\frac{S_n - \mathsf{E}(S_n)}{\sqrt{\mathsf{D}(S_n)}} \implies {\cal N}(0,1).
\end{equation}
\end{en-theorem}

\begin{en-corollary}\label{Cor2}
If  $\sup_{n} C_n = C < \infty$ and  
\begin{equation}\label{Db>0}
\inf_{i,n}
\mathsf{D}(f^{(n)}_i(X^{(n)}_i)) \geq c > 0,
\end{equation}
and the conditions  (\ref{ijm0})  along with 
\begin{equation}\label{cor2-e1}
\lim_{n \to \infty} n \alpha_n (\alpha^\beta_n)^2 = \infty,
\end{equation}
hold, then convergence  (\ref{CLT}) is valid. 
\end{en-corollary}
Conditions (\ref{Db>0}), (\ref{ijm0}), and (\ref{cor2-e1}) replace classical Dobrushin's conditions~(\ref{Dobr-cond}). 

\fi

\ifen
\section{Auxiliary results}\label{sec4}
The proofs of the propositions and lemmata here are similar to those in \cite{Hall-Heyde} and \cite{Seth-Varadhan} with certain necessary changes due to the new setting. We highlight that in the meanwhile the result is, indeed, new. Also, the lemmata \ref{lem1}, \ref{lem2}, and \ref{lem4} are new. The notations are mostly similar to those in  \cite{Seth-Varadhan}. 
In particular, denote by  $\|Z\|_{L^{\infty}}$ the essential supremum norm of the random variable:
$$
\|Z\|_{L^{\infty}} := \mathop{\text{ess}}\limits_{\mathsf P}\sup_{\omega\in \Omega} |Z(\omega)|.
$$

\begin{en-proposition}[\cite{Seth-Varadhan}, Proposition 3.1, with reference to \cite{Hall-Heyde}]\label{pro2}
Let for any $n \geq 1$ the process $\{(W^{(n)}_i: 0 \leq i \leq n \}$ be a martingale with respect to a filtration $({\cal G}^{(n)}_i)$,  $W^{(n)}_0 = 0$, and let 
$\xi^{(n)}_i := W^{(n)}_i - W^{(n)}_{i-1}$. If 
\begin{equation}\label{clt_a}
\max_{1 \leq i \leq n} \|\xi^{(n)}_i\|_{L^\infty} \to 0 , \quad \sum_{i = 1}^n \mathsf{E}[(\xi^{(n)}_i)^2 | G^{(n)}_i] \to 1 \quad \text{in} \quad L^2
\end{equation}
then weak convergence holds, 
$$
W^{(n)}_n \implies {\cal N}(0,1).
$$
\end{en-proposition}

More general martigale CLT results see in 
\cite[Theorem 3.2]{Hall-Heyde} and
\cite[Theorem 5.5.8]{Lip-Shi}, and most general in \cite{Jac-Shi}. 

\medskip

Everywhere in what follows, as in  \cite{Seth-Varadhan}  it will be assumed that the functions  $\{ f^{(n)}_i\}$ are such that 
$$
\mathsf{E}[f^{(n)}_i(X^{(n)}_i)] = 0, \quad 1 \leq i \leq n, \; n \geq 1.
$$
Let
\begin{equation*}
Z^{(n)}_k := \sum_{i=k}^n \mathsf{E}[f^{(n)}_i(X^{(n)}_i)| X^{(n)}_k],
\end{equation*}
so that 
\begin{align*}
Z^{(n)}_k = 
\begin{cases}
f^{(n)}_k(X^{(n)}_k) + \sum_{i=k+1}^n \mathsf{E}[f^{(n)}_i(X^{(n)}_{i})| X^{(n)}_k],  \quad 1 \leq k \leq n-1, 
  \\  \\ 
f^{(n)}_n (X^{(n)}_n), \quad  k = n.
\end{cases}
\end{align*}
Then, for any $1 \leq k \leq n-1$, the representation holds, 

\begin{equation*}
f^{(n)}_k(X^{(n)}_k) = Z^{(n)}_k - \mathsf{E}[Z^{(n)}_{k+1}| X^{(n)}_k].  
\end{equation*}
More than that, for any $2 \leq k \leq n-1$ the same expression may be equivalently represented as  $(Z^{(n)}_k - \mathsf{E}[Z^{(n)}_{k}| X^{(n)}_{k-1}]) + (\mathsf{E}[Z^{(n)}_{k}| X^{(n)}_{k-1}] - \mathsf{E}[Z^{(n)}_{k+1}| X^{(n)}_k])$. It follows that the sum $S_n$ admits the following, crucial for what follows martingale representation, 
$$
S_n = \sum_{i=1}^n  f^{(n)}_i(X^{(n)}_i) = Z^{(n)}_1 + \sum_{k=2}^n [Z^{(n)}_k - \mathsf{E}[Z^{(n)}_{k}| X^{(n)}_{k-1}]].
$$
This is a transformation proposed by M.I. Gordin in \cite{Gordin}.  
Since all summands in the right-hand side of the latter formula are non-correlated, we have, 
$$
\mathsf{D}(S_n) = \sum_{k=2}^n \mathsf{D}(Z^{(n)}_k - \mathsf{E}[Z^{(n)}_{k}| X^{(n)}_{k-1}]) + \mathsf{D}(Z^{(n)}_1).
$$ 
Let
\begin{equation}\label{XiEq}
\xi^{(n)}_k = \frac{1}{\mathsf{D}(S_n)}[Z^{(n)}_k - \mathsf{E}[Z^{(n)}_{k}| X^{(n)}_{k-1}]].
\end{equation}
Then the process $M^{(n)}_k = \sum_{\ell=1}^{k}  \xi^{(n)}_\ell$ is a martingale with respect to the filtration $\mathcal{F}^{(n)}_k = \sigma\{X^{(n)}_\ell :1 \leq \ell \leq k\}$ as $1 \leq k \leq n$. So, the task is to verify both conditions in  (\ref{clt_a}). 

Further, by  $f^{(n)}_i\pi_{i,j}f^{(n)}_j$ we denote the product of the functions  $f^{(n)}_i$ and $\pi_{i,j}f^{(n)}_j$. 

\medskip

Notice that for any $i<j$, $n$,  and function $f^{(n)}_j$, 
\begin{equation*}
\operatorname{Osc}(\pi_{i,j}(f^{(n)}_j)) \le \delta(\pi_{i,j})\operatorname{Osc}(f^{(n)}_j).
\end{equation*}

\begin{en-lemma}\label{lem1}
Under the conditions of the theorem, for any  $1 \leq i \leq j \leq n $, 
\begin{equation*}
\|\pi_{i,j} f^{(n)}_j\|_{B} \leq 2C_n (1 - \alpha^\beta_n)^{\kappa^\beta_{j} - \kappa^\beta_{i}}, \quad \operatorname{Osc}(\pi_{i,j}(f^{(n)}_j)^2) \leq 2 C^2_n (1 - \alpha^\beta_n)^{\kappa^\beta_{j} - \kappa^\beta_{i}}; 
\end{equation*}
and also for  $1 \leq l < i \leq j \leq n $,

$$
\operatorname{Osc}(\pi_{l,i}(f^{(n)}_i\pi_{i,j}f^{(n)}_j)) \leq 6 C^2_n (1 - \alpha^{\beta}_n)^{\kappa^\beta_{i} - \kappa^\beta_{l}}(1 - \alpha^{\beta}_n)^{\kappa^\beta_{j} - \kappa^\beta_{i}};
$$
    
\end{en-lemma}

The next result from \cite{Seth-Varadhan} provides a lower bound for the variance of the sum.
Formally, our assumptions differ from those in \cite{Seth-Varadhan}; however, the proof of this sole proposition is just based on the same common conditions.

\begin{en-proposition}[\cite{Seth-Varadhan}, proposition 3.2]\label{{pro3}}
Under the assumptions of the theorem \ref{thm1}, 
for any  $n \geq 1$,
$$
\mathsf{D}(S_n) \geq \frac{\alpha_n}{4} \sum_{i=1}^n  \mathsf{D}(f^{(n)}_i(X^{(n)}_i)).
$$
\end{en-proposition}

To prove this proposition, the following two auxiliary statements are helpful. We state them here for the completeness of the presentation.

\begin{en-lemma}[\cite{Seth-Varadhan}, lemma 4.1]\label{lem02}
Let $\lambda$ be a probability measure on $\mathbf{X} \times \mathbf{X}$ with marginals $\nu_1$ and $\nu_2$ respectively. Let $\pi\left(x_{1}, d x_{2}\right)$ and $\widehat{\pi}\left(x_{2}, d x_{1}\right)$ be the corresponding transition probabilities in the two directions so that $\nu_1 \pi=\nu_2$ and $\nu_2 \hat{\pi}=\nu_1$. Let $f\left(x_{1}\right)$ and $g\left(x_{2}\right)$ be square integrable with respect to $\nu_1$ and $\nu_2$, respectively. If
$$
\int f\left(x_{1}\right) \nu_1\left(d x_{1}\right)=\int g\left(x_{2}\right) \nu_2\left(d x_{2}\right)=0
$$
then,
$$
\left|\int f\left(x_{1}\right) g\left(x_{2}\right) \lambda\left(d x_{1}, d x_{2}\right)\right| \leq \sqrt{\delta(\pi)}\|f\|_{L_{2}(\nu_1)}\|g\|_{L_{2}(\nu_2)} .
$$
\end{en-lemma}

\begin{en-lemma}[\cite{Seth-Varadhan}, lemma 4.2]\label{lem03}
 Let $f\left(x_{1}\right)$ and $g\left(x_{2}\right)$ be square integrable with respect to $\nu_1$ and $\nu_2$, respectively. Then,
$$
E\left[\left(f\left(x_{1}\right)-g\left(x_{2}\right)\right)^{2}\right] \geq \alpha(\pi) \mathsf{D}\left(f\left(x_{1}\right)\right)
$$
as well as
$$
E\left[\left(f\left(x_{1}\right)-g\left(x_{2}\right)\right)^{2}\right] \geq\alpha(\pi) \mathsf{D}\left(g\left(x_{2}\right)\right) .
$$
\end{en-lemma}

Further, for the ratio $\|Z^{(n)}_k\|_{L^{\infty}}/\sqrt{\mathsf{D}(S_n)}$, we have an estimate in the following lemma, which is the analogue of \cite[lemma 3.2]{Seth-Varadhan}.

\begin{en-lemma}
\label{lem2}
Under the assumptions of the theorem \ref{thm1}, the equality holds, 
\begin{equation*}
\lim_{n \to \infty} \sup_{1 \leq k \leq n} \frac{\|Z^{(n)}_k\|_{L^{\infty}}}{\sqrt{\mathsf{D}(S_n)}} = 0.
\end{equation*}
\end{en-lemma}

\begin{en-lemma}[\cite{Seth-Varadhan}, lemma 3.3] \label{lem3}
Let $\left\{Y_l^{(n)}: 1 \leq l \leq n\right\}$ and $\left\{\mathcal{G}_l^{(n)}: 1 \leq l \leq n\right\}$, for $n \geq 1$, are, respectively, the sequence of non-negative random variables and sigma-fields such that  $\sigma\left\{Y_1^{(n)}, \ldots, Y_l^{(n)}\right\} \subset \mathcal{G}_l^{(n)}$. Suppose that 
$$
\lim _{n \rightarrow \infty} \mathsf{E}\left[\sum_{l=1}^n Y_l^{(n)}\right]=1 \quad \text {and} \quad \sup _{1 \leq i \leq n}\left\|Y_i^{(n)}\right\|_{L^{\infty}} \leq \epsilon_n,
$$
where $\lim _{n \rightarrow \infty} \epsilon_n=0$. Suppose in addition that 
$$
\lim _{n \rightarrow \infty} \sup _{1 \leq l \leq n-1} \operatorname{Osc}\left(\mathsf{E}\left[\sum_{j=l+1}^n Y_j^{(n)} \mid \mathcal{G}_l^{(n)}\right]\right)=0.
$$
Then
$$
\lim _{n \rightarrow \infty} \sum_{l=1}^n Y_l^{(n)}=1 \quad \text \quad \mbox{in} \quad L^2.
$$
\end{en-lemma}
{\em (To verify it! Are the conditions sufficient for the convergence in $L_2$?)}

Further, denote 
$$
v_j^{(n)} :=  \mathsf{E}[(\xi_j^{(n)})^2 \mid \mathcal{F}_{j-1}^{(n)}],
$$ 
where $\xi_j^{(n)}$ are defined in (\ref{XiEq}).
These random variables are measurable with respect to the sigma-fields  $\mathcal{G}^{(n)}_j = \mathcal{F}^{(n)}_{j-1}$ for $2 \leq j \leq n $, respectively. The last lemma is the analogue of \cite[lemma 3.4]{Seth-Varadhan}

\begin{en-lemma} \label{lem4}
Under the assumptions of the theorem \ref{thm1} the convergence 
$$
\sup _{2 \leq l \leq n-1}\, 
\operatorname{Osc}\left(\mathsf{E}\left[\sum_{j=l+1}^n v_j^{(n)} \mid \mathcal{F}_{l-1}^{(n)}\right](\omega)\right)=o(1), \quad n\to\infty,
$$
holds true.
\end{en-lemma}
{\em NB: Here 
the expression $o(1)$ is non-random according to the definition in   (\ref{Osc-rv}).}

\fi

\ifen

\section{Proof of the theorem  1}\label{sec5}

{\em Proof.}
It follows from the lemma \ref{lem2} that it suffices to show convergence 
$$
M_n^{(n)} / \sqrt{\mathsf D \left(S_n\right)} \implies {\cal N}(0,1).
$$
This property follows from the proposition  \ref{pro2}, because it was already verified that 
\begin{enumerate}
\item[(a)] $\sup _{2 \leq k \leq n}\left\|\xi_k^{(n)}\right\|_{L^{\infty}} \rightarrow 0, \quad n\to\infty$;
    
\item[(b)]  $\mathsf E \left\vert \sum_{k=2}^n \mathsf{E}\left[\left(\xi_k^{(n)}\right)^2 \mid \mathcal{F}_{k-1}^{(n)}\right] - 1\right\vert^2 \rightarrow 0, \quad
n\to \infty$.
       
\end{enumerate}
Here, part (a) follows from the lemma \ref{lem2}, while part (b) follows from the lemma \ref{lem3} and from part (a). The theorem is proved. 
\hfill $\square$

\fi

\ifen

\fi

\ifen
\section*{Financial support}
For both authors this work was supported by the Foundation for the Advancement of Theoretical Physics and Mathematics ``BASIS''.

\fi

\ifru

\fi

\ifen

\fi


\end{document}